\documentclass[a4paper,11pt]{extarticle}

\usepackage[left=3cm,right=1.5cm,top=2cm,bottom=2cm,bindingoffset=0cm]{geometry}
\usepackage{setspace} % двойной интервал

\usepackage{cmap}					% поиск в PDF
\usepackage[T2A]{fontenc}			% кодировка
\usepackage[utf8]{inputenc}			% кодировка исходного текста
\usepackage[english]{babel}	% локализация и переносы

\usepackage{geometry} % пакет для установки полей
\usepackage{mathrsfs} % Math font
\usepackage{graphicx}
\usepackage{amsmath,amsthm,amsfonts,amssymb} %AMS
\usepackage{icomma} % Умная запятая 

\usepackage{array,tabularx,tabulary,booktabs} % Дополнительная работа с таблицами
\usepackage{longtable}  % Длинные таблицы
\usepackage{multirow} % Слияние строк в таблице
\usepackage{hhline}

\usepackage{graphicx}  % Для вставки рисунков
\graphicspath{{images/}}  % папки с картинками
\usepackage{wrapfig} % Обтекание рисунков и таблиц текстом
\usepackage{epstopdf} %eps в pdf

\usepackage[usenames]{color}
\usepackage{color}
\usepackage{colortbl}

\usepackage{titlesec} %добавление точки после номера главы
\titleformat{\section}{\filcenter\normalfont\Large\bfseries}{\thesection.}{0.2em}{} %добавление точки после номера главы и центрирование название главы

\usepackage{indentfirst}

\usepackage{cite}

\usepackage{hyperref}

\usepackage{algorithm}
\usepackage{algpseudocode}

\begin{document}

\begin{center}
%\textbf{ИДЕНТИФИЦИРУЕМОСТЬ ПРОСТРАНСТВЕННОЙ SEIR-HCD МОДЕЛИ РАСПРОСТРАНЕНИЯ COVID-19 \footnote[1]{Работа выполнена в рамках государственного задания Министерства науки и высшего образования РФ (тема "Аналитическое и численное исследование обратных задач об определении параметров источников атмосферного или водного загрязнения и (или) параметров среды", код темы: FENG-2023-0004) и при поддержке Российского научного фонда (проект № 23-71-10068).}}
\textbf{ IDENTIFIABILITY OF THE SPATIAL SEIR-HCD MODEL OF COVID-19 PROPAGATION\footnote{The research was carried out within the state assignment of Ministry of Science and Higher Education of the Russian Federation (theme №~FENG-2023-0004, ``Analytical and numerical study of inverse problems on recovering parameters of atmosphere or water pollution sources and (or) parameters of media'').}}

%\textbf{О.~И.~Криворотько, Т.~А.~Звонарева, А.~В.~Неверов}

\vspace{1em}
\textbf{Olga Krivorotko$^{1,2a}$, Tatiana Zvonareva$^{2b}$, Andrei Neverov$^{2c}$}

%\textit{ $^1$Югорский государственный университет, \\
%ул. Чехова, 16, г. Ханты-Мансийск 628012, Россия,\\
%$^2$Институт математики им.~С.~Л.~Соболева СО РАН, \\
%пр-кт. Академика Коптюга, 4, г. Новосибирск 630090, Россия\\
%e-mail: $^a$krivorotko.olya@mail.ru, $^b$t.zvonareva@g.nsu.ru, $^c$a.neverov@g.nsu.ru}

\vspace{1em}
\textit{
$^1$Yugra University, Khanty-Mansiysk, 628012 Russia,\\
$^2$Sobolev Institute of Mathematics SB RAS, Novosibirsk, 630090 Russia\\
e-mail: $^a$krivorotko.olya@mail.ru, $^b$t.zvonareva@g.nsu.ru, $^c$a.neverov@g.nsu.ru}
\end{center}

%В работе исследуется вопрос идентифицируемости пространственной математической модели распространения быстротекущих эпидемий, основанной на законе действующих масс и диффузионных процессах. Алгоритм исследования основан на глобальных методах анализа чувствительности Соболя и байесовском подходе, которые в совокупности позволили уменьшить границы вариации неизвестных параметров для дальнейшего решения задачи идентификации параметров по измерениям количества выявленных случаев заболевания, критических и умерших. Показано, что для идентификации диффузионных коэффициентов, отвечающих за скорость перемещения индивидуумов в пространстве, необходимо использовать дополнительную информацию о процессе.
This paper investigates the identifiability of a spatial mathematical model of the spread of fast-moving epidemics based on the law of acting masses and diffusion processes. The research algorithm is based on global methods of Sobol sensitivity analysis and Bayesian approach, which together allowed to reduce the variation boundaries of unknown parameters for further solving the problem of parameter identification by measurements of the number of detected cases, critical and dead. It is shown that for identification of diffusion coefficients responsible for the rate of movement of individuals in space, it is necessary to use additional information about the process.

%\textbf{Ключевые слова:} модель <<реакции-диффузии>>, идентифицируемость, анализ чувствительности, задача об источнике, оптимизация.

\vspace{1em}

\section*{INTRODUCTION}

%Математическое моделирование распространения инфекционных заболеваний далеко продвинулось в результате вспышки новой коронавирусной инфекции, вызванной вирусом SARS-CoV-2 в 2019 году в г. Ухань. Начиная с классических дифференциальных моделей эпидемиологии, основа которых была заложена в работе Кермака и МакКендрика в 1927 году~\cite{Kermack_McKendrick_1927}, и заканчивая агентно-ориентированными моделями и моделями управления средним полем, а также их комбинацией, исследователи до сих пор сталкиваются с неполными и неточными данными для устойчивого описания распространения эпидемий. Более полный обзор моделей распространения быстротекущих эпидемий приведен в работе~\cite{KOI_KSI_2024} и кратко на диаграмме~\ref{fig:Model_scheme}.

\begin{figure}[!h]
    \centering
    \includegraphics[width=1\textwidth]{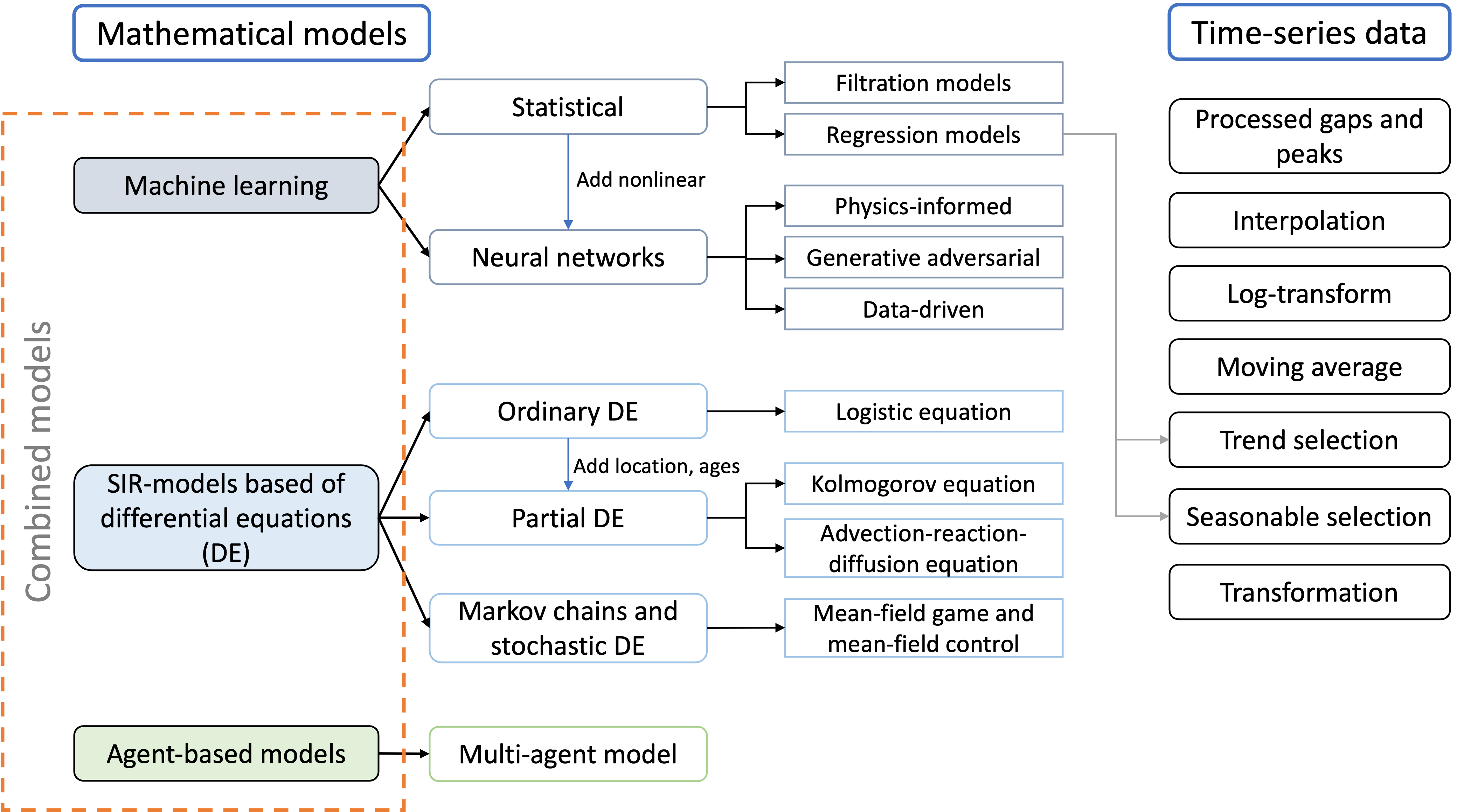}
    \caption{Brief review of mathematical models of epidemic spread and relationship to statistical approaches.}
    \label{fig:Model_scheme}
\end{figure}

Mathematical modelling of the spread of infectious diseases has come a long way with the outbreak of a new coronavirus infection caused by SARS-CoV-2 virus in 2019 in Wuhan. From classical differential models of epidemiology, whose foundation was laid in the work of Kermack and McKendrick in 1927~\cite{Kermack_McKendrick_1927}, to agent-based models and mean-field control models, as well as their combination, researchers are still faced with incomplete and inaccurate data to sustainably describe the spread of epidemics. A more comprehensive review of models of fast-moving epidemic spread is given in the paper~\cite{KOI_KSI_2024} and summarized in the figure~\ref{fig:Model_scheme}.

%\caption{Краткий обзор математических моделей распространения эпидемий и взаимосвязь со подходами обработки статистических данных.}

%В данной работе исследуются прямая и обратная задача для пространственной SEIR-HCD модели распространения COVID-19 в регионе~\cite{KOI_KSI_2020}:
%\begin{eqnarray}\label{eq:model_vector}
%\dfrac{\partial u}{\partial t} = \nabla (n v \nabla u) + f(u, q),\quad u(x,0)=u_0(x), \qquad t\in (0, T).
%\end{eqnarray}
%Здесь $u(x,t) = (s, e, i, r, h, c, d)(x,t)$~-- вектор состояний системы, характеризующий плотности восприимчивых ($s$), бессимптомных инфицированных ($e$), больных COVID-19 ($i$), выздоровевших ($r$), госпитализированных ($h$), критических случаев, требующих подключение аппарата искусственной вентиляции легких (ИВЛ) ($c$) и умерших ($d$) в результате COVID-19. Параметры $v=(v_s, v_e, v_i, v_r)$ и $q(x)\in C(\mathbb{R}^m)$ являются непрерывными пространственными функциями и описывают скорость перемещения соответствующих групп населения в пространстве и особенности распространения инфекции в популяции, соответственно. Без ограничения общности будем рассматривать безразмерное одномерное пространство по $x$, то есть $x\in (0,1)$.

This paper investigates the direct and inverse problem for a spatial SEIR-HCD model of COVID-19 propagation in the region\cite{KOI_KSI_2020}:
\begin{eqnarray}\label{eq:model_vector}
\dfrac{\partial u}{\partial t} = \nabla (n v \nabla u) + f(u, q),\quad u(x,0)=u_0(x), \quad t\in (0, T).
\end{eqnarray}
Here $u(x,t) = (s, e, i, r, h, c, d)(x,t)$ is the vector of system states characterizing the densities of susceptible ($s$), asymptomatic infected ($e$), COVID-19 patients ($i$), recovered ($r$), hospitalized ($h$), critical cases requiring ALV device connection ($c$) and died ($d$) as a result of COVID-19. The parameters $v=(v_s, v_e, v_i, v_r)$ and $q(x)\in C(\mathbb{R}^m)$ are continuous spatial functions that describe the rate of movement of the respective populations in space and the patterns of infection in the population, respectively. Without limiting generality, we will consider a dimensionless one-dimensional space on $x$, i.~e., $x\in (0,1)$.

%Математическая модель~(\ref{eq:model_vector}) основывается на законе баланса масс, при котором вся моделируемая популяция по времени постоянна, т.е. $n(x,t) = n(x)$, где $n(x,t) = s(x,t) + e(x,t) + i(x,t) + r(x,t) + h(x,t) + c(x,t) + d(x,t)$.

The mathematical model~(\ref{eq:model_vector}) is based on a mass balance law in which the entire modelled population is constant in time, i.~e. $n(x,t) = n(x)$, where 
\[n(x,t) = s(x,t) + e(x,t) + i(x,t) + r(x,t) + h(x,t) + c(x,t) + d(x,t).\]

%Впервые в 1937 году А.Н.~Колмогоров, И.Г.~Петровский и Н.С.~Пискунов~\cite{KolmPetrPiskunov_1937} предложили и обосновали математическую модель~(\ref{eq:model_vector}), получившую впоследствии название модели <<реакции-диффузии>>, а также строго доказали, что если начальное условие удовлетворяет следующим ограничениям

%то при $f(u)=u(1-u)$ динамика изменения популяции во времени описывается скоростью $v^*=2\sqrt{f^\prime(0)nv}$. Интересно отметить, что схожие волновые решения обнаружены в распределенной модели распространения эпидемий~\cite{Murray_2007}. В указанной ситуации такие колебательные решения интерпретируются как повторные вспышки эпидемий.

For the first time in 1937, A.~N.~Kolmogorov, I.~G.~Petrovsky and N.~S.~Piskunov~\cite{KolmPetrPiskunov_1937} proposed and justified the mathematical model~(\ref{eq:model_vector}), which later became known as the reaction-diffusion model, and also strictly proved that if the initial condition satisfies the following restrictions
\[0\leqslant u(x,0)\leqslant 1,\quad u(x,0)=0\;\;\; \forall\, x<x_1 \quad \text{and}\quad u(x,0)=1\;\;\; \forall\, x>x_2\geqslant x_1,\]
then for $f(u)=u(1-u)$ the population dynamics in time is described by the rate $v^*=2\sqrt{f^\prime(0)nv}$. It is interesting to note that similar wave solutions are found in the distributed epidemic propagation model~\cite{Murray_2007}. In the above situation, such oscillatory solutions are interpreted as repeated outbreaks of epidemics.

Taking spatial heterogeneity into account makes it possible to more accurately model an epidemic from a spreading center (a large city in a country, a capital city in a region, etc.) under known initial conditions. Thus, in~\cite{Viguerie_2020, Aristov_2021, Barwolff_2021, Lau_2022}, estimates of the spread of COVID-19 in the first months of the epidemic, taking into account passenger traffic, were obtained. In~\cite{Boscheri_2021}, the authors propose a multiscale spatial model based on a system of SEIR-type kinetic transport equations describing the passenger population moving at large scales (outside the city) and a system of diffusion equations describing the urban population acting at small scales. This modelling approach avoids unrealistic effects of traditional diffusion models in epidemiology, such as infinite diffusion rates at large scales and mass migration dynamics. The system parameters include both epidemiological characteristics (number of contacts, disease rate) and socioeconomic characteristics (social distance, vaccination). Given empirically defined parameters, the authors obtained the spatial distribution of the COVID-19 epidemic in Italy by solving systems of partial derivative equations using the finite volume method on unstructured meshes. In~\cite{Takacs_2021}, an integro-differential equation with delay describing the epidemic spreading process is proposed and studied, in which the solution has some biologically valid features. The feature of epidemiological models consists of the principle of transmission from an infected $I$ to a susceptible $S$ individual, which is characterized by the type of function $f(u, q)$, i.~e. the function
\[f_I(S,I,q) = \dfrac{\beta(x) S I^p}{1+\kappa I}\]
characterizes the non-linear transmission of viral infection with a number of contacts $\beta(x)$ and a parameter $\kappa>0$ describing social distance or other socioeconomic constraints during the period of infection spread.

%Статья организована следующим образом. В разделе~\ref{sec_dir_pr} приведена постановка SEIR-HCD модели, ее пространственная интерпретация и методы ее численного решения. В разделе~\ref{sec_inv_pr} приведена постановка обратной задачи для пространственной SEIR-HCD модели. Анализ чувствительности двумя глобальными методами проведен в разделе~\ref{sec_identifiability}. Алгоритм численного решения обратной задачи с учетом результатов идентифицируемости описан в разделе~\ref{sec_algorithm}.

The paper is organized as follows. In section~\ref{sec_dir_pr}, the formulation of the SEIR-HCD model, its spatial interpretation, and methods for its numerical solution are introduced. In section~\ref{sec_inv_pr}, the formulation of the inverse problem for the spatial SEIR-HCD model is presented. Sensitivity analyses by two global methods are carried out in section~\ref{sec_identifiability}. The algorithm for the numerical solution of the inverse problem, taking into account the identifiability results, is described in~\ref{sec_algorithm}.

\section{DIRECT PROBLEM FORMULATION FOR THE SPATIAL SEIR-HCD MODEL}\label{sec_dir_pr}

%В работе~\cite{KOI_ZNY_Eurasian_2022} сформулирована и проанализирована математическая модель распространения вспышки COVID-19, основанная на законе баланса масс и описываемая системой из 7 обыкновенных дифференциальных уравнений:

In~\cite{KOI_ZNY_Eurasian_2022}, a mathematical model of the COVID-19 outbreak propagation, based on the mass balance law and described by a system of 7 ordinary differential equations, is formulated and analyzed:
\begin{eqnarray}\label{eq:SEIR-HCD_ODE}
    \left\{\begin{array}{ll}
        \dfrac{dS}{dt} = -\dfrac{\alpha_I S(t)I(t)}{N} - \dfrac{\alpha_E S(t)E(t)}{N} + \dfrac{1}{t_{imm}}R(t), \\
        \dfrac{dE}{dt} = \dfrac{\alpha_I S(t)I(t)}{N} + \dfrac{\alpha_E S(t)E(t)}{N} - \dfrac{1}{t_{inc}}E(t), \\
        \dfrac{dI}{dt} = \dfrac{1}{t_{inc}}E(t) - \dfrac{1}{t_{inf}}I(t), \\
        \dfrac{dR}{dt} = \dfrac{\beta(t)}{t_{inf}}I(t) + \dfrac{1-\varepsilon_{HC}}{t_{hosp}}H(t) - \dfrac{1}{t_{imm}}R(t), \\
        \dfrac{dH}{dt} = \dfrac{1-\beta(t)}{t_{inf}}I(t) + \dfrac{1-\mu}{t_{crit}}C(t) - \dfrac{1}{t_{hosp}}H(t), \\
        \dfrac{dC}{dt} =  \dfrac{\varepsilon_{HC}}{t_{hosp}}H(t) - \dfrac{1}{t_{crit}}C(t), \\
        \dfrac{dD}{dt} =  \dfrac{\mu}{t_{crit}}C(t).
    \end{array}\right.
\end{eqnarray}

%В модели~(\ref{eq:SEIR-HCD_ODE}) не учитывается индекс самоизоляции от Яндекса, так как с 2021 года он перестал обновляться ввиду снятия ограничительных и контрольных мер. В модели учитываются:
%\begin{itemize}
%    \item инкубационный период $t_{inc}$ бессимптомного течения инфекции,
%    \item возможность повторного заражения ввиду ослабления иммунитета в течение времени $t_{imm}$,
%    \item регламентированные в моделируемом регионе время госпитализации $t_{hosp}$ и длительность использования аппарата ИВЛ $t_{crit}$,
%    \item возможная смертность после пребывания в критическом состоянии.
%\end{itemize}

The model~(\ref{eq:SEIR-HCD_ODE}) does not take into account the self-isolation index from Yandex, as it has stopped being updated since 2021 due to the removal of restrictive and control measures. The model takes into account:
\begin{itemize}
    \item incubation period $t_{inc}$ of asymptomatic infection,
    \item the possibility of re-infection due to weakening of immunity during the time $t_{imm}$,
    \item hospitalization time $t_{hosp}$ and duration of artificial lung ventilation (ALV) device use $t_{crit}$ regulated in the simulated region,
    \item possible mortality after critical state.
\end{itemize}

%Описание параметров модели~(\ref{eq:SEIR-HCD_ODE}) и начальных условий, характеризующих распространение COVID-19 в Новосибирской области с 31.01.2022 (период вспышки эпидемии ввиду появления нового штамма <<Омикрон>> в регионе), приведено в Таблице~\ref{tab_parameters}.

The description of parameters of the model~(\ref{eq:SEIR-HCD_ODE}) and initial conditions characterizing the spread of COVID-19 in the Novosibirsk region from 31.01.2022 (the period of the outbreak due to the appearance of a new Omicron strain in the region) is given in the table~\ref{tab_parameters}.

\begin{table}[h!]
\centering
\caption{SEIR-HCD model parameters characterizing the distribution of COVID-19 in the Novosibirsk region from 31.01.2022.}\label{tab_parameters}
 \begin{tabular}{|l|p{9.8cm}|l|} 
 \hline
 Parameter & Description & Value \\ 
  \hline
$\alpha_E$ & Infection parameter for asymptomatic and susceptible groups & 0.0922 \\
$\alpha_I$ & Infection parameter for infected and susceptible groups & 0.3856 \\
$\beta(t)$ & Proportion of individuals with late IgG antibodies to SARS-CoV-2 obtained by Invitro & Data~\cite{covid_data}\\
$\varepsilon_{HC}$ & Proportion of hospitalized cases with ALV support & 0.0376 \\
$\mu$ & Proportion of dead from COVID-19 & 0.4754 \\
$t_{inc}$ & Incubation period duration (days) & $\mathrm{LogN} (4.6,4.8)$ \\
$t_{inf}$ & Infection period duration (days) & $\mathrm{LogN} (6.6, 4.9)$ \\
$t_{hosp}$ & Hospitalization period duration (days) & $\mathrm{LogN}(3, 7.4)$ \\
$t_{crit}$ & ALV support duration (days) & $\mathrm{LogN}(6.2, 1.7)$ \\
$t_{imm}$ & Average duration of natural immunity (days) & $\mathrm{LogN}(150, 30)$ \\
$N$ & Population of the Novosibirsk region & 2798170 \\
$S_0$ & Initial number of susceptible cases & 2734917 \\
$E_0$ & Initial number of asymptomatic cases & 4329 \\
$I_0$ & Initial number of infected cases & 3508 \\
$R_0$ & Initial number of recovered cases & 32333 \\
$H_0$ & Initial number of hospitalized cases & 219 \\
$C_0$ & Initial number of critically ill cases & 54 \\
$D_0$ & Initial number of dead cases & 4932 \\
\hline
\end{tabular}
\end{table}

%Пространственная SEIR-HCD модель~\cite{Barwolff_2021} распространения вспышки COVID-19 основывается на модели~(\ref{eq:SEIR-HCD_ODE}) и описывается системой уравнений типа "реакции-диффузии":

%Здесь $v = (v_s, v_e, v_i, v_r)$~--- скорости перемещения в пространстве групп восприимчивых, бессимптомных и инфицированных носителей COVID-19, а также выздоровевших от COVID-19, и $q = (\alpha_i, \alpha_e, \varepsilon_{HC}, \mu)$~--- вектор параметров модели, характеризующие распространение COVID-19 в конкретном регионе. В работе предполагается, что госпитализированная группа населения, критическая и умершие не могут перемещаться в пространстве.

The spatial SEIR-HCD model~\cite{Barwolff_2021} of the COVID-19 outbreak propagation is based on the model~(\ref{eq:SEIR-HCD_ODE}) and is described by a system of reaction-diffusion type equations
\begin{equation}\label{eq:PDE_SEIR-HCD}
\left\{ \begin{array}{l}
    \partial_t s = \nabla (n\,v_s\nabla s) - \alpha_i si - \alpha_e se + t_{imm}^{-1}r, \\[3pt]
    \partial_t e = \nabla (n\,v_e\nabla e) + \alpha_i si + \alpha_e se - t_{inc}^{-1}e, \\[3pt]
    \partial_t i = \nabla (n\,v_i\nabla i) + t_{inc}^{-1}e - t_{inf}^{-1} i, \\[3pt]
    \partial_t r = \nabla (n\,v_r\nabla r) + \beta(t) t_{inf}^{-1}i + (1-\varepsilon_{HC})t_{hosp}^{-1}h - t_{imm}^{-1}r, \\[3pt]
    \partial_t h = (1-\beta(t))t_{inf}^{-1}i + (1-\mu)t_{crit}^{-1}c - t_{hosp}^{-1}h, \\[3pt]
    \partial_t c = \varepsilon_{HC}t_{hosp}^{-1}h - t_{crit}^{-1}c, \\[3pt]
    \partial_t d = \mu t_{crit}^{-1} c,
\end{array}\right.
\end{equation}
initial conditions
\begin{equation}\label{init_cond}
\begin{array}{l}
    s(x,\,0) = e^{-(x + 1)^4} + e^{-\frac{(x - 0.35)^2}{10^{-2}}} + \dfrac{1}{8}\left(e^{-\frac{(x - 0.62)^4}{10^{-5}}} + e^{-\frac{(x - 0.52)^4}{10^{-5}}} + e^{-\frac{(x - 0.42)^4}{10^{-5}}}\right) + \dfrac{1}{4}e^{-\frac{(x - 0.735)^4}{10^{-5}}}, \\
    e(x,\,0) = \dfrac{1}{20}e^{-\frac{(x - 0.75)^4}{10^{-5}}}, \quad i(x,\,0) = i_0, \quad r(x,\,0) = r_0,\quad h(x,\,0) = h_0,\\
    c(x,\,0) = c_0,\quad d(x,\,0) = d_0, 
\end{array}
\end{equation}
where $u_0 = \dfrac{U_0}{N}$, and boundary conditions of the form
\begin{equation}\label{bound_cond}
\begin{array}{l}
    \partial_x u(0,\,t) = 0, \quad u(1,\,t) = 0.
\end{array}
\end{equation}
Here $v = (v_s, v_e, v_i, v_r)$ are the velocities of the groups of susceptible, asymptomatic and infected COVID-19 carriers, and recovered of COVID-19, and $q = (\alpha_i, \alpha_e, \varepsilon_{HC}, \mu)$ is the vector of model parameters characterizing the spread of COVID-19 in a particular region. In this paper, it is assumed that the hospitalized population, critical and dead populations cannot move in space.

%Под прямой задачей для модели~(\ref{eq:PDE_SEIR-HCD})--(\ref{bound_cond}) в данной работе понимается задача моделирования процесса распространения COVID-19, в которой требуется найти векторную функцию $u(x,t)$ при заданных $v$ и $q$.

The direct problem for the model~(\ref{eq:PDE_SEIR-HCD})--(\ref{bound_cond}) in this paper denotes the problem of modelling the COVID-19 propagation process, which requires finding the vector function $u(x,t)$ when $v$ and $q$ are given.

%Параметрами модели~(\ref{eq:PDE_SEIR-HCD})--(\ref{bound_cond}) является пара $(v, q)$, где $v$~--- вектор скалярных параметров диффузии с единицами измерения 1/чел.день.

The parameters of the model~(\ref{eq:PDE_SEIR-HCD})--(\ref{bound_cond}) are the pair $(v, q)$, where $v$ is the vector of scalar diffusion parameters with units of 1/(person per day).

\subsection{Numerical methods for solving the direct problem}

%В разделе приведены численные методы решения прямой начально-краевой задачи~(\ref{eq:PDE_SEIR-HCD})--(\ref{bound_cond}): методы конечных элементов (МКЭ) и конечных разностей (МКР).

In this section, we present numerical methods for solving the direct initial boundary value problem~(\ref{eq:PDE_SEIR-HCD})--(\ref{bound_cond}): finite element method (FEM) and finite difference method (FDM).

\subsubsection{Finite element method}

%Для МКЭ расчетная область разбивается на подобласти, называемые конечными элементами, внутри которых функция $u(x,t)$ приближается выбранными базисными функциями.

For FEM, the computational domain is partitioned into sub-regions called finite elements, within which the function $u(x,t)$ is approximated by selected basis functions.

%Сначала рассматривается слабая постановка задачи, для этого исходное уравнение умножается на пробную функцию $\phi$ и интегрируется по пространству. Рассмотрим на примере первого уравнения из системы (\ref{eq:PDE_SEIR-HCD}):

First we consider a weak formulation of the problem, for this purpose the initial equation is multiplied by the trial function $\phi$ and integrated over the space. Let us consider the first equation from the system~(\ref{eq:PDE_SEIR-HCD}) as an example:
\[\int_0^1\psi\partial_{t}sdx=\int_0^1\psi\left(\nabla\cdot(n\,v_{s}\nabla s)-\beta_{i}si-\beta_{e}se\right)dx.\]

%Интегрируя по частям слагаемое с градиентом, получаем слабую постановку задачи:

Integrating by parts the summand with the gradient, we obtain a weak formulation of the problem:
\begin{equation}\label{eq:FEMweak}
  \int_0^1\psi\partial_{t}sdx=\int_0^1\left(-\nabla\psi(n\,v_{s}\nabla s)-\psi\beta_{i}si-\psi\beta_{e}se\right)dx + \psi n\,v_{s}\nabla s|_{x=0}^1.
\end{equation}

%В качестве базисных функций на элементе выбраны линейные функции, тогда в соответствие к ним выбраны пробные функции вида:

%где $x_i < x_{i+1}, i=0..N_{x}$ - границы конечных элементов. Тогда решение можно представить в виде:

%где $u_k = (s_k,e_k,i_k,r_k,h_k,c_k,d_k)$ - вектор значений функций в точке $x_k$, момент времени $t$. Подставляя функции в таком виде в уравнение (\ref{eq:FEMweak}) и заменяя производную по времени конечно-разностным аналогом, получаем:

%откуда получаем значения локальной матрицы $K$ жесткости, которые должны быть такие, что

Linear functions are chosen as basis functions on the element, then trial functions are chosen in correspondence to them:
\[ \psi_{k}=\begin{cases}
\frac{x-x_{k-1}}{x_{k}-x_{k-1}}, & x\in[x_{k-1,}x_{k}],\\
\frac{x_{k+1}-x}{x_{k+1}-x_{k}}, & x\in[x_{k,}x_{k+1}],\\
0 & ,\text{ otherwise.}
\end{cases}\]
where $x_i < x_{i+1}$, $i=0,\dots,N_{x}$ are the boundaries of the finite elements. Then the solution can be represented as
\[\sum_{k=0}^{N_x} u_k \psi_k,\]
where $u_k = (s_k,e_k,i_k,r_k,h_k,c_k,d_k)$ is the vector of function values at the point $x_k$ at time $t$. Substituting the functions in this form into equation~(\ref{eq:FEMweak}) and replacing the time derivative with a finite-difference analogue, we obtain
\begin{multline*}
    \int_{0}^{1}\psi_{s}\psi_{k}\frac{s_k-s_{k-1}}{\tau}dx=\sum_{k=1}^{N_{x}-1}\int_{0}^{1}\left(\nabla\psi_{k}\nabla\psi_{s}(n_k\,v_{s}s_k)-\psi_{s}\psi_{k}\beta_{i}s_k i_k-\psi_{s}{\psi_{k}\beta}_{e}s_k e_k\right)dx \\
+(\psi_{s}\psi_{N}n_N\,v_{s}\nabla s_N)|_{x=1}-(\psi_{s}\psi_{0}n_0\,v_{s}\nabla s_0)|_{x=0},
\end{multline*}
whence we obtain the values of the local stiffness matrix $K$, which must be such that
\begin{multline*}
[K_j u = -n_{k}\,v_{s}\frac{1}{h}\left(\begin{array}{cc}
1 & -1\\
-1 & 1
\end{array}\right)\left(\begin{array}{c}
s_{j}\\
s_{j+1}
\end{array}\right)-\frac{1}{\tau}\frac{h}{6}\left(\begin{array}{cc}
2 & 1\\
1 & 2
\end{array}\right)\left(\begin{array}{c}
s_{j}\\
s_{j+1}
\end{array}\right)\\
+\beta_{i}s_{j}\frac{h}{6}\left(\begin{array}{cc}
2 & 1\\
1 & 2
\end{array}\right)\left(\begin{array}{c}
i_{j}\\
i_{j+1}
\end{array}\right)+\beta_{e}s_{j}\frac{h}{6}\left(\begin{array}{cc}
2 & 1\\
1 & 2
\end{array}\right)\left(\begin{array}{c}
e_{j}\\
e_{j+1}
\end{array}\right).
\end{multline*}

%И правая часть рассчитывается как:

And the right side is calculated as
\[M_{j}u=-\frac{1}{\tau}\frac{h}{6}\left(\begin{array}{cc}
2 & 1\\
1 & 2
\end{array}\right)\left(\begin{array}{c}
s_{j}^{t-\tau}\\
s_{j+1}^{t-\tau}
\end{array}\right).\]

%Складывая все матрицы $K_j$ и правые части, получаем систему линейных алгебраических уравнений, решением которой будут искомые коэффициенты $u_k$.

Adding all matrices $K_j$ and right-hand sides, we obtain a system of linear algebraic equations, the solution of which will be the required coefficients $u_k$.

\subsubsection{Finite difference method}

%МКР приближает непрерывный вектор функций $u(x,t)$ его сеточным аналогом $(u_k^j)$. Для этого вводится сетка в замкнутой области $\Omega = \{ (x,t) \; | \; 0\le x\le 1, \; 0\le t\le T \}$:

The FDM approximates a continuous vector of functions $u(x,t)$ by its grid analogue $(u_k^j)$. For this purpose, a grid in the closed region $\Omega = \{ (x,t) \; | \; 0\le x\le 1, \; 0\le t\le T \}$ is introduced:
\begin{align*}
    \omega = \{ \left( x_{k}, t_{j} \right) \; | \; x_{k} = kh, \; t_{j} = j\tau, \; j=0,\dots,N_{x}, \; k=0,\dots,N_{t} \}, 
\end{align*}
where $h = \dfrac{1}{N_{x}}$ and $\tau = \dfrac{T}{N_{t}}$.

%Далее первые и вторые производные приближаются конечными разностями с порядком аппроксимации $O(\tau+h^2)$, в результате чего получаются разностные уравнения:

%для $k = 1,\dots,N_x-1$, $j = 0,\dots,N_t-1$. А разностные аналоги начальных и граничных условий имеют вид:

Then the first and second derivatives are approximated by finite differences with approximation order $O(\tau+h^2)$, resulting in difference equations
\begin{equation*}
\begin{array}{rl}
    \dfrac{s_k^{j + 1} - s_k^j}{\tau} = & \negthickspace\negthickspace v_s\dfrac{n_{k + 1}^j - n_{k - 1}^j}{2h_x}\dfrac{s_{k + 1}^j - s_{k - 1}^j}{2h_x} + v_s n_k^j\dfrac{s_{k + 1}^j - 2s_k^j + s_{k - 1}^j}{h_x^2} - \alpha_i s_k^j i_k^j - \alpha_e s_k^j e_k^j + \\
    & \negthickspace\negthickspace +\; t_{imm}^{-1}r_k^j, \\
    \dfrac{e_k^{j + 1} - e_k^j}{\tau} = & \negthickspace\negthickspace v_e\dfrac{n_{k + 1}^j - n_{k - 1}^j}{2h_x}\dfrac{e_{k + 1}^j - e_{k - 1}^j}{2h_x} + v_e n_k^j\dfrac{e_{k + 1}^j - 2e_k^j + e_{k - 1}^j}{h_x^2} + \alpha_i s_k^j i_k^j + \alpha_e s_k^j e_k^j - \\
    & \negthickspace\negthickspace -\; t_{inc}^{-1}e_k^j, \\
    \dfrac{i_k^{j + 1} - i_k^j}{\tau} = & \negthickspace\negthickspace v_i\dfrac{n_{k + 1}^j - n_{k - 1}^j}{2h_x}\dfrac{i_{k + 1}^j - i_{k - 1}^j}{2h_x} + v_i n_k^j\dfrac{i_{k + 1}^j - 2i_k^j + i_{k - 1}^j}{h_x^2} + t_{inc}^{-1}e_k^j - t_{inf}^{-1} i_k^j, \\
    \dfrac{r_k^{j + 1} - r_k^j}{\tau} = & \negthickspace\negthickspace v_r\dfrac{n_{k + 1}^j - n_{k - 1}^j}{2h_x}\dfrac{r_{k + 1}^j - r_{k - 1}^j}{2h_x} + v_r n_k^j\dfrac{r_{k + 1}^j - 2r_k^j + r_{k - 1}^j}{h_x^2} + \beta(t) t_{inf}^{-1}i_k^j + \\
    & \negthickspace\negthickspace +\; (1 - \varepsilon_{HC})t_{hosp}^{-1}h_k^j - t_{imm}^{-1}r_k^j, \\
    \dfrac{h_k^{j + 1} - h_k^j}{\tau} = & \negthickspace\negthickspace (1 - \beta(t))t_{inf}^{-1}i_k^j + (1 - \mu)t_{crit}^{-1}c_k^j - t_{hosp}^{-1}h_k^j, \\
    \dfrac{c_k^{j + 1} - c_k^j}{\tau} = & \negthickspace\negthickspace  \varepsilon_{HC}t_{hosp}^{-1}h_k^j - t_{crit}^{-1}c_k^j, \\
    \dfrac{d_k^{j + 1} - d_k^j}{\tau} = & \negthickspace\negthickspace \mu t_{crit}^{-1} c_k^j,
\end{array}
\end{equation*}
for $k = 1,\dots,N_x-1$, $j = 0,\dots,N_t-1$, and the difference analogues of the initial and boundary conditions are of the form:
\begin{equation*}
\begin{array}{l}
    s_k^0 = s(x_k,0), \quad e_k^0 = e(x_k,0), \quad i_k^0 = i_0, \\ r_k^0 = r_0, \quad h_k^0 = h_0, \quad c_k^0 = c_0, \quad d_k^0 = d_0, 
\end{array}
    \quad \text{for} \quad k = 0,\dots,N_x,
\end{equation*}
\begin{equation*}
\begin{array}{l}
    \dfrac{-3u_{0}^{j + 1} + 4u_{1}^{j + 1} - u_{2}^{j + 1}}{2h} = 0, \text{ for } j = 0,\dots,N_t-1, \\
    u_{N_x}^j = 0, \text{ for } j = 1,\dots,N_t.
\end{array}
\end{equation*}

%Откуда получаются явные выражения для определения функции $u_k^j$ для всех $k$ и $j$. 

Whence we obtain explicit expressions for the definition of the function $u_k^j$ for all $k$ and $j$.

\section{INVERSE PROBLEM STATEMENT}\label{sec_inv_pr}

%В обратной задаче для модели~(\ref{eq:PDE_SEIR-HCD})-(\ref{bound_cond}) помимо вектора функций $u(x,t)$ неизвестными являются начальные функции $s(x,0)$, $e(x,0)$ и $i(x,0)$.

In the inverse problem for the model~(\ref{eq:PDE_SEIR-HCD})--(\ref{bound_cond}), in addition to the function vector $u(x,t)$, the unknowns are the initial functions $s(x,0)$, $e(x,0)$, and $i(x,0)$.

%Полагаем, что задана дополнительная информация о процессе при $K$ днях измерений и $x_i$ точках потенциальных источников, т.е.

%Здесь $I_{k}$, $C_{k}$ и $D_{k}$~-- количество выявленных, критических и умерших случаев в день $k$ соответственно.

We assume that additional information about the process at $K$ days of measurements and $x_i$ points of potential sources is given, i.~e.
\begin{align}\label{ip_data_PDE_SEIR-HCD}
\begin{array}{c}
    I(t_{k}) = I_{k}, \quad C(t_{k}) = C_{k}, \quad D(t_{k}) = D_{k}, \quad k=1,\dots,K.
\end{array}
\end{align}
Here $I_{k}$, $C_{k}$, and $D_{k}$ are the number of infected, critical, and dead cases on day $k$, respectively.

%Задача определения вектора трех функций в общем виде по информации типа~(\ref{ip_data_PDE_SEIR-HCD}) является некорректной (т.е. ее решение может быть неединственным и/или неустойчивым)~\cite{Kabanikhin_2009}. Поэтому предполагается, что количество источников, имеющие экспоненциальный вид~(\ref{init_cond}) (аналог <<гауссовой шапочки>>), известно (ими могут быть крупные центры в регионах, города в стране) и имеет вид~\cite{Viguerie_2020}:

The problem of determining the vector of three functions in general form from information of the type~(\ref{ip_data_PDE_SEIR-HCD}) is incorrect (i.~e., its solution may be non-unique and/or unstable)~\cite{Kabanikhin_2009}. Therefore, it is assumed that the number of sources having exponential form~(\ref{init_cond}) (analogue of ``Gaussian cap'') is known (they can be large centers in regions, cities in the country) and has the form~\cite{Viguerie_2020}
\begin{gather*}
    s(x,0) = a_1^s e^{-\frac{(x-b_1^s)^4}{c_1^s}} + a_2^s e^{-\frac{(x-b_2^s)^4}{c_2^s}} + a_3^s e^{-\frac{(x-b_3^s)^4}{c_3^s}},\\
    e(x,0) = a_1^e e^{-\frac{(x-b_1^e)^4}{c_1^e}} + a_2^e e^{-\frac{(x-b_2^e)^4}{c_2^e}} + a_3^e e^{-\frac{(x-b_3^e)^4}{c_3^e}},\\
    i(x,0) = i_0.
\end{gather*}

%Обратная задача~(\ref{eq:PDE_SEIR-HCD})-(\ref{bound_cond}),~(\ref{ip_data_PDE_SEIR-HCD}) состоит в восстановлении вектора параметров $q = (a^p,b^p,c^p,i_0)$, $p = \{s, e\}$, модели~(\ref{eq:PDE_SEIR-HCD}) по данным $I_{k}$, $C_{k}$ и $D_{k}$ вида~(\ref{ip_data_PDE_SEIR-HCD}).

The inverse problem~(\ref{eq:PDE_SEIR-HCD})-(\ref{bound_cond}),~(\ref{ip_data_PDE_SEIR-HCD}) consists of recovering the parameter vector $q = (a^p,b^p,c^p,i_0)$, $p = \{s, e\}$, model~(\ref{eq:PDE_SEIR-HCD}) using data $I_{k}$, $C_{k}$, and $D_{k}$ of the form~(\ref{ip_data_PDE_SEIR-HCD}).

\subsection{Variational formulation of the problem} 

%Для уточнения вектора неизвестных параметров $q$ формулируется вариационная постановка обратной задачи, состоящая в минимизации квадратичного целевого функционала:

To refine the vector of unknown parameters $q$, we consider a variational formulation of the inverse problem, which consists of minimizing the quadratic target functional
\begin{equation}\label{eq:IP_func}
    J(q) = \sum_{k=1}^{K} |I(t_k;q) - I_{k}|^2 + |C(t_k;q) - C_{k}|^2 + |D(t_k;q) - D_{k}|^2.
\end{equation}

\subsection{Sensitivity analysis of the inverse problem parameters}\label{sec_identifiability}

\subsubsection{Sobol sensitivity analysis}

%Один из способов исследования чувствительности модели к неизвестным параметрам является метод на основе анализа дисперсии модели \cite{Cukier, Sobol}. В основе данных методов лежит рассмотрение распределения неизвестных параметров с помощью методов Монте-Карло, главной задачей которого является подсчет индексов чувствительности для каждого из исследуемых параметров. Пусть модель задана в виде:

%где $q_1,..., q_k$ -- набор неизвестных параметров. Сгенерируем матрицу \textbf{$Q$} размерности $N \times (k+2)$, которая является случайным набором неизвестных параметров $\vec q = {q_i, i=1...k}$ в заданных неуточненных границах. Тогда основанный на дисперсии эффект первого порядка на модель (\ref{eq1}) для параметра $q_i$ можно записать как:

%Здесь $q_i$ -- $i$-й параметр, $\textbf{Q}_{~i}$ -- сгенерированная матрица неизвестных параметров без $q_i$. Смысл оператора математического ожидания $E_{\textbf{Q}_{~i}}(Y|q_i)$ заключается в том, что среднее значение для $Y$ берется по всем возможным значениям  $\textbf{Q}_{~i}$ при фиксированном $q_i$. В то же время дисперсия $V_{q_i}$ берется по всем возможным значениям $q_i$. Соответствующая мера чувствительности параметра $q_i$ (коэффициент чувствительности первого порядка) записывается следующим образом:

%Здесь $V(Y)$ - дисперсия по всем строкам матрицы $Y$. Таким образом, индексы чувствительности представляют собой ожидаемое сокращение дисперсии, которое было бы получено, если бы значение параметра $q_i$ можно было зафиксировать, нормированное на общую дисперсию.

One of the ways to investigate the sensitivity of a model to unknown parameters is the method based on analysis of variance of the model~\cite{Cukier, Sobol}. These methods are based on the consideration of the distribution of unknown parameters using Monte Carlo methods, the main task of which is to calculate sensitivity indices for each of the parameters under study. Let the model be given in the form
\begin{equation}\label{eq1}
Y = f(q_1,..., q_k),
\end{equation}
where $q_1,..., q_k$ is the set of unknown parameters. We generate a matrix \textbf{$Q$} of dimension $N \times (k+2)$, which is a random set of unknown parameters $\vec q = {q_i, i=1...k}$ in the given unspecified bounds. Then the variance-based first-order effect on the model~(\ref{eq1}) for parameter $q_i$ can be written as
\begin{equation}
    V_{q_i}(E_{\textbf{Q}_{~i}}(Y|q_i)).
\end{equation}
Here $q_i$ is the $i$th parameter, $\textbf{Q}_{~i}$ is the generated matrix of unknown parameters without $q_i$. The meaning of the expectation operator $E_{\textbf{Q}_{~i}}(Y|q_i)$ is that the mean value for $Y$ is taken over all possible values of $\textbf{Q}_{~i}$ at fixed $q_i$. At the same time, the variance of $V_{q_i}$ is taken over all possible values of $q_i$. The corresponding sensitivity measure of parameter $q_i$ (first-order sensitivity coefficient) is written as follows:
\begin{equation}
    S_i = \frac{V_{q_i}(E_{\textbf{Q}_{~i}}(Y|q_i))}{V(Y)}.
\end{equation}
Here $V(Y)$ is the variance across all rows of the matrix $Y$. Thus, the sensitivity indices represent the expected reduction in variance that would be obtained if the value of parameter $q_i$ could be fixed, normalised by the total variance.

%Подробно с программной реализацией и исходным кодом можно ознакомиться в документации \url{https://salib.readthedocs.io/en/latest/basics.html}, а также в работе \cite{Saltelli_2010}. 

Details of the software implementation and source code can be found in the documentation \url{https://salib.readthedocs.io/en/latest/basics.html} and in the paper~\cite{Saltelli_2010}. 

%Для модели (\ref{eq:PDE_SEIR-HCD}) был проведен анализ идентифицируемости неизвестных параметров $\vec q = (\alpha_i, \alpha_e, t_{inc}, t_{inf}, \beta(t), \epsilon_{HC}, t_{hosp}, t_{imm}, \mu, t_{crit}, v_s, v_e, v_i, v_r)$ с помощью метода Соболя в реализации библиотеки Salib на языке Python \url{https://github.com/SALib/SALib}.

For the model~(\ref{eq:PDE_SEIR-HCD}), an identifiability analysis of the unknown parameters $\vec q = (\alpha_i, \alpha_e, t_{inc}, t_{inf}$, $\beta(t), \varepsilon_{HC}, t_{hosp}, t_{imm}, \mu, t_{crit}, v_s, v_e, v_i, v_r)$ using the Sobol method in the Python implementation of the Salib library \url{https://github.com/SALib/SALib}.

%На рисунке~\ref{Sobol_timeline} представлены значения индексов чувствительности $S_i$ для всех неизвестных параметров в различные временные срезы, т.е. исследовалась дисперсия системы в пяти временных точках решения прямой задачи: на 40-ой, 80-ый, 120-ый, 160-ый и 200-ый дни моделирования. 

\begin{figure}[h!]
    \centering
    \includegraphics[width=1\textwidth]{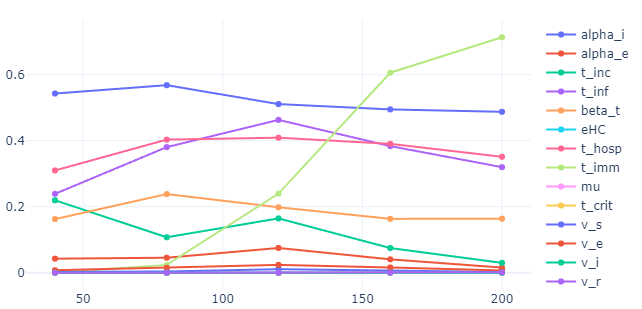}
    \caption{Sensitivity indices for unknown parameters $\vec q$ over time. Values are presented for the 40th, 80th, 120th, 160th and 200th days of modelling.}
    \label{Sobol_timeline}
\end{figure}

Figure~\ref{Sobol_timeline} shows the values of the sensitivity indices $S_i$ for all unknown parameters at different time slices, i.~e., the variance of the system at five time points of the direct problem solution was investigated: on the 40th, 80th, 120th, 160th, and 200th days of the simulation. 

%По графику можно сделать вывод, что параметр $t_{imm}$ со временем становится более влиятельным для системы, что естественно, поскольку популяционный иммунитет возрастает со временем. В то же время индексы чувствительности параметров $\alpha_i, t_{inc}, t_{inf}, t_{hosp}$ наоборот постепенно уменьшаются. Это может быть вызвано специфичностью модели -- камеры, отвечающие за распространение заболевания ($I, E, H$), в уравнениях которых участвуют последние параметры, имеют волнообразную особенность и со временем затухают (см. рисунок \ref{Forward_problem}). 

\begin{figure}[h!]
    \centering
    \includegraphics[width=1\textwidth]{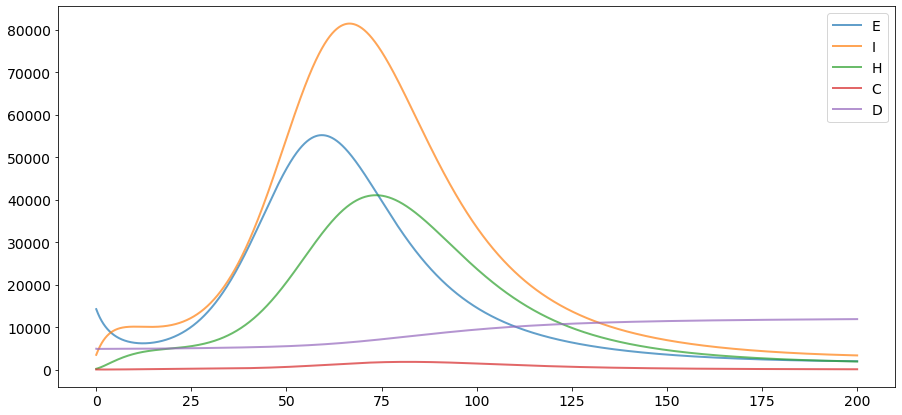}
    \caption{Solution of the direct problem for the model~(\ref{eq:PDE_SEIR-HCD}) with given parameters $\alpha_i = 0.3856$, $\alpha_e = 0.0922$, $t_{inc} = 5$,  $t_{inf} = 8$, $\beta_t = 0.4$, $\epsilon_{HC} = 0.0376$, $t_{hosp} = 7$, $t_{imm} = 175$, $\mu = 0.4754$, $t_{crit} = 9$, $v_s = 5e-5$, $v_e = 1e-3$, $v_i = 1e-10$, $v_r = 5e-5$.}
    \label{Forward_problem}
\end{figure}

The graph shows that the parameter $t_{imm}$ becomes more influential for the system with time, which is natural because population immunity increases with time. At the same time, the sensitivity indices of the parameters $\alpha_i, t_{inc}, t_{inf}, t_{hosp}$ on the contrary gradually decrease. This may be due to the specificity of the model: the compartments responsible for the spread of the disease ($I, E, H$), in whose equations the latter parameters are involved, have a wave-like feature and decay with time (see figure~\ref{Forward_problem}). 

%На рисунке (\ref{Sobol_days}) представлены столбчатые диаграммы значений индексов чувствительности в зависимости от времени. По диаграмме можно сделать заключение о том, что $\alpha_i$ является самым чувствительным параметром, в то время как параметры $v_s, v_i, v_r, v_e$ являются наименее чувствительными и не привносят существенного вклада на дисперсию результатов модели (\ref{eq:PDE_SEIR-HCD}).

\begin{figure}[h!]
\begin{minipage}[h]{0.47\linewidth}
\center{\includegraphics[width=1\linewidth]{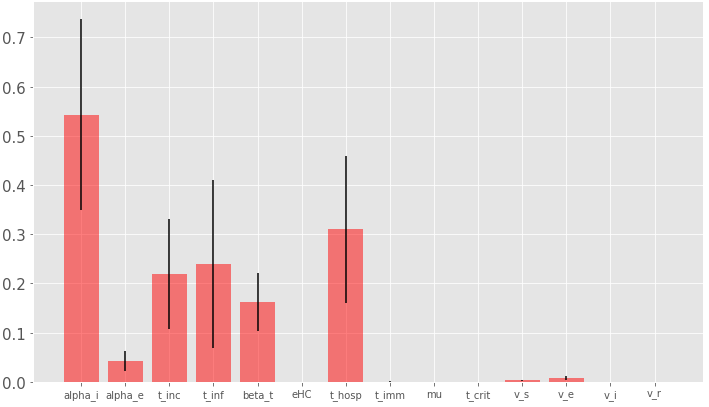}} a) $T=40$ days \\
\end{minipage}
\hfill
\begin{minipage}[h]{0.47\linewidth}
\center{\includegraphics[width=1\linewidth]{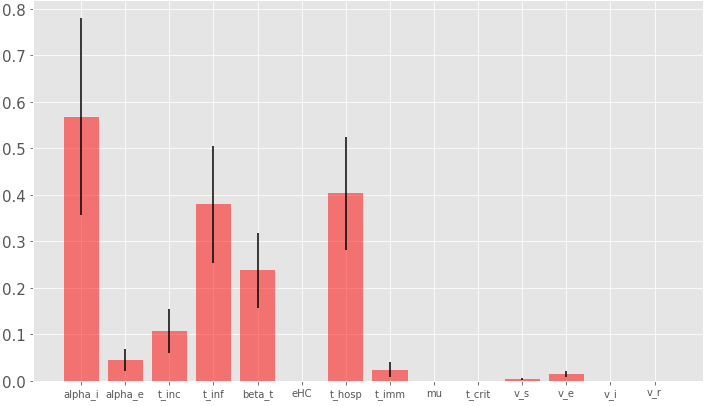}} \\b) $T=80$ days
\end{minipage}
\vfill
\begin{minipage}[h]{0.47\linewidth}
\center{\includegraphics[width=1\linewidth]{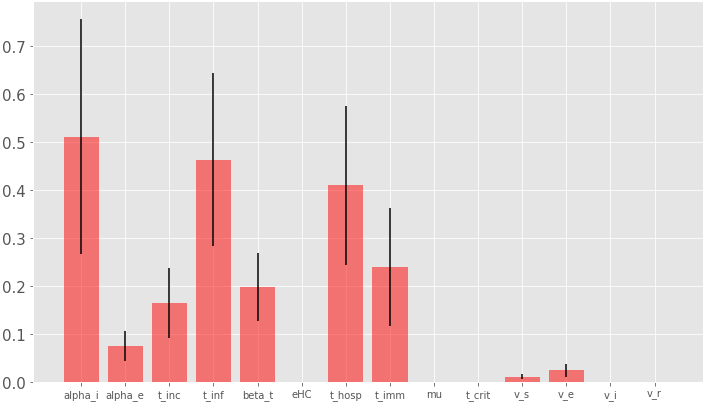}} c) $T=120$ days \\
\end{minipage}
\hfill
\begin{minipage}[h]{0.47\linewidth}
\center{\includegraphics[width=1\linewidth]{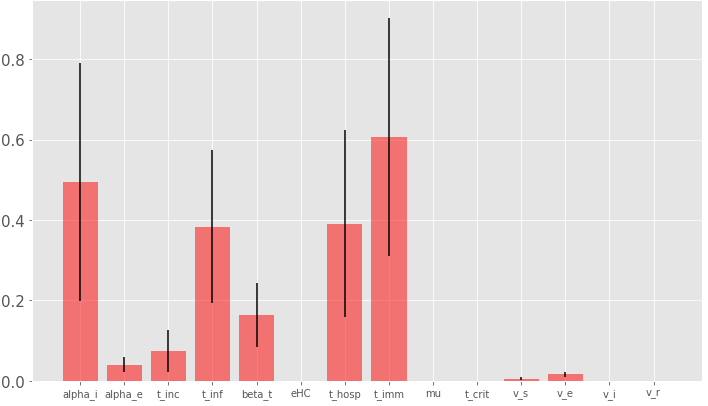}} d) $T=160$ days \\
\end{minipage}
\caption{Diagrams of sensitivity index values for unknown parameters $\vec q$ for different time slices. The red bars represent the values of indices $S_i$, the black lines represent the confidence intervals for $S_i$.}
\label{Sobol_days}
\end{figure}

The figure~(\ref{Sobol_days}) shows bar charts of the values of sensitivity indices as a function of time. The diagram shows that $\alpha_i$ is the most sensitive parameter, while the parameters $v_s, v_i, v_r, v_e$ are the least sensitive and do not contribute significantly to the variance of the results of the model~(\ref{eq:PDE_SEIR-HCD}).

\subsubsection{Method of sensitivity analysis based on Bayesian approach}

%Иной метод исследования чувствительности неизвестных параметров к реальным данным и уточнения их границ перед поиском оптимальных значений в ходе решения обратной задачи является метод на основе байесовского подхода~\cite{emulator}. Его основная идея состоит в построении простой гибридной модели (эмулятора) для исследуемой  модели (симулятора). В силу особенности вида эмулятора, он требует гораздо меньше вычислительных и временных ресурсов, поэтому может быть запущен многократно для оценки степени влияния параметров на результат модели.

Another method of investigating the sensitivity of unknown parameters to real data and refining their bounds before searching for optimal values in the course of solving the inverse problem is the method based on the Bayesian approach~\cite{emulator}. Its main idea is to construct a simple hybrid model (emulator) for the model under study (simulator). Due to the feature of the emulator type, it requires much less computational and time resources, so it can be run repeatedly to assess the degree of influence of parameters on the model result.

%Эмулятор представляет собой следующую функцию от неизвестных параметров:

%где $h_i(\vec q)$ -- мономы порядка $i$, $\beta_i$ -- регрессионные коэффициенты, $u(\vec q)$ -- гауссовский процесс.

The emulator is the following function of the unknown parameters:
\begin{equation}
        g(\vec q) = \sum_{i=1}^{p} h_i(\vec q)\cdot\beta_i+u(\vec q),
\end{equation}
where $h_i(\vec q)$ are monomials of order $i$, $\beta_i$ are regression coefficients, $u(\vec q)$ is a Gaussian process.

The general algorithm for investigating the identifiability of parameters using the emulator is as follows:
\begin{enumerate}
    \item We specify a set of 250 parameters $\vec q = \{q_i | q_i\in [a_i, b_i], i=1...N$\} (here $N$ is the number of unknown parameters, $[a_i, b_i]$ is the unspecified bounds for them) using the Latin hypercube~\cite{LHD} for a uniform distribution of points in $N$-dimensional space.
    \item The direct problem~(\ref{eq:PDE_SEIR-HCD}) is solved 250 times from the unknown parameter values generated in step 1 and the model results are stored.
    \item An emulator with zero mean function and exponential correlation function of Gaussian process is created:
    \begin{equation}
        c(\vec q,\vec q' ) = \sigma^2\cdot\exp \left[ - \sum_{i=1}^{p}\frac{(q_i-q'_i)^2}{\delta_i^2} \right].
    \end{equation}
    \item We implement the selection of $\beta_i$ parameters of the emulator to match the simulator output data that were obtained in (2). For this purpose, we solve the inverse problem using the quasi-Newtonian L-BFGS-B~\cite{lfbgs} optimization method.
    \item The last step in the analysis is the determination of the parameter indentifiability space, which is carried out using the history matching method. To do this, 50 000 parameter sets are generated and then the emulator outputs at these points are compared to real observations to determine data-driven plausibility bounds. History matching involves the calculation of an implausibility metric, which establishes how well a particular set of input parameters describes real observations:
    \begin{equation}\label{metric}
        I(\vec{q^*}) = \frac{|z - E(g(\vec{q^*}))|}{\sqrt{Var[z - E(g(\vec{q^*}))]}}.
    \end{equation}
    Here $z$ is the real observations (statistical data), $\vec{q^*}$ is the parameter set, $E(g(\vec{q^*}))$ is the mean of the emulator.
    \item After computing the implausibility metric, the parameter space is divided into two parts based on the threshold value $t$ ($t=3$ by the rule $3\sigma$~\cite{Pukelsheim}). Thus, we are interested in points with $I(\vec q) < 3$, which form the boundaries of the plausible space.
    \item In the last step of the algorithm, the different bounds for all outputs are crossed.
\end{enumerate}

%Было произведено исследование чувствительности параметров к реальным данным, а также уточнение границ с помощью построения эмулятора для 14 неизвестных параметров модели (\ref{eq:PDE_SEIR-HCD}):

The sensitivity of the parameters to real data was investigated, and the bounds were refined by constructing an emulator for 14 unknown model parameters~(\ref{eq:PDE_SEIR-HCD})
\begin{equation}
     \vec q = (\alpha_i, \alpha_e, t_{inc}, t_{inf}, \beta(t), \epsilon_{HC}, t_{hosp}, t_{imm}, \mu, t_{crit}, v_s, v_e, v_i, v_r).
\end{equation}
 
%На вход алгоритма были переданы следующие неуточненные границы параметров, в которых происходила симуляция точек методом Латинского гиперкуба (таблица~\ref{boundaries}).

The following unspecified parameter bounds were passed to the input of the algorithm, in which point simulation using the Latin Hypercube method took place (table~\ref{boundaries}).

 \begin{table}[h!]
\begin{center}
\caption{Unspecified bounds on the parameters of the model~\ref{eq:PDE_SEIR-HCD} in which sensitivity analyses were performed using emulator construction.}
\label{boundaries}
\begin{tabular}{c|c} 
\hline
Parameter & Unspecified boundaries \\ 
\hline
$\alpha_i$ & [0.0; 38.9] \\
$\alpha_e$ &  [0.0; 9.3] \\
$t_{inc}$ & [0.0; 505.0]  \\
$t_{inf}$ & [0.0; 808.0] \\
$\beta_t$ & [0.0; 40.4] \\
$\epsilon_{HC}$ & [0.0; 3.8] \\
$t_{hosp}$ & [0.0; 707.0] \\
$t_{imm}$ & [0.0; 17675.0] \\
$\mu$ & [0.0; 48.0] \\
$t_{crit}$ & [0.0; 909.0] \\
$v_s$ & [0.0; 0.005] \\
$v_e$ & [0.0; 0.101] \\
$v_i$ & [0.0; 0.001] \\
$v_r$ & [0.0; 0.005] \\
\hline
\end{tabular}
\end{center}
\end{table}

%В качестве измеряемых данных для сравнения выходов эмуляторы использовались реальные данные по Новосибирской области по количеству госпитализированных, выздоровевших и умерших от COVID-19. Таким образом было построено 3 пространства правдоподобности по алгоритму выше для каждой из статистик, которые были пересечены между собой.

Real data for Novosibirsk region on the number of hospitalized, recovered and died from COVID-19 were used as measured data to compare the outputs of the emulators. Thus, 3 plausibility spaces were constructed using the algorithm above for each of the statistics, which were overlapped with each other.
 
%На рисунке \ref{emulator_nroy} представлен результат исследования чувствительности неизвестных параметров к реальным данным (для более точного сравнения исследуемые границы были отнормированы на интервал [0,1]). 

\begin{figure}[H]
    \centering
    \includegraphics[width=\textwidth]{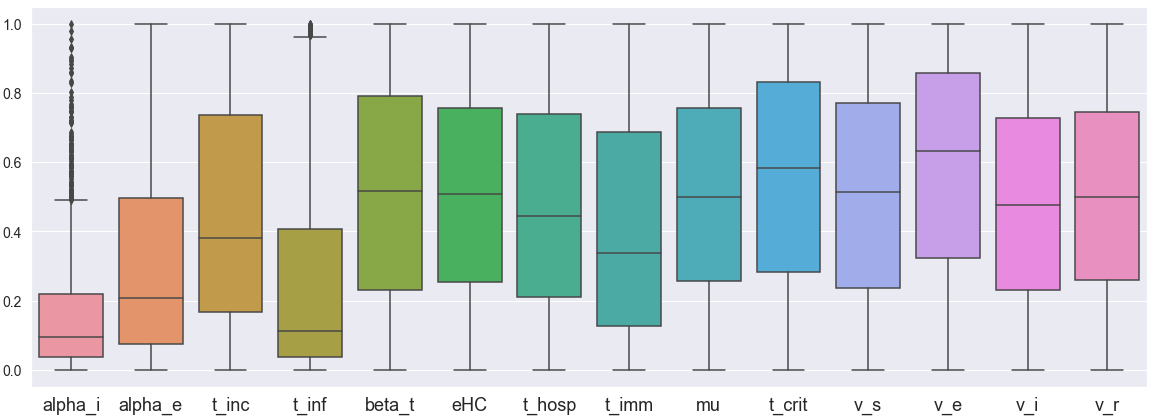}
    \caption{Diagrams of density distributions of unknown parameters in the plausibility space obtained using the Bayesian approach. The boundaries of ``boxes with whiskers'' represent 25, 50 and 75 quantiles of distributions, the lines at the bottom and top are the unspecified boundaries of unknown parameters.}
\label{emulator_nroy}
\end{figure}

Figure~\ref{emulator_nroy} shows the result of the study of the sensitivity of unknown parameters to real data (for a more accurate comparison, the investigated boundaries were normalized to the interval [0,1]).
 
%Можно заметить, что с помощью данного подхода в качестве наиболее чувствительных параметров были выделены $\alpha_i, \alpha_e, t_{inf}, t_{crit}, t_{imm}$, т.к. их 50-е квантили смещены к одной из заданных границ, соответственно точки в пространстве правдоподобности густо сконцентрированы в смещенной области. В свою очередь остальные параметры являются менее чувствительными, поскольку для них точки в правдоподобном пространстве распределены равномерно, и их  уточненные границы не могут быть специфично определены по реальным данным.

It can be seen that using this approach, $\alpha_i, \alpha_e, t_{inf}, t_{crit}, t_{imm}$ were selected as the most sensitive parameters, since their 50th quantiles are shifted to one of the specified boundaries, and, accordingly, the points in the plausibility space are densely concentrated in the shifted region. In turn, the other parameters are less sensitive, because for them the points in a plausible space are uniformly distributed, and their refined bounds cannot be specifically determined from real data.

\section{ALGORITHM FOR SOLVING INVERSE PROBLEM}\label{sec_algorithm}

%Задачу минимизации целевого функционала планируется решать методом глобальной оптимизации тензорного поезда (Tensor Train global optimization, TT)~\cite{Zheltkova_2018}, алгоритм которого представлен далее.
%\begin{algorithm}
%\renewcommand{\thealgorithm}{\hspace{-5pt}}
%\caption{метода ТТ}
%\begin{algorithmic}[1]
%    \Require{Нижняя и верхняя границы пространства решений $b_{min}$ и $b_{max}$, число параметров (размерность пространства решений) $d$, количество узлов по всем направлениям $n$, максимально возможный ранг вагонов $r_{max}$, количество итераций $N_{TT}$, начальный сдвиг функционала $\alpha$, функция отображения $h(J(q) - \alpha)$.}
%    \State Ввести сетку с $n$ узлами по каждому из $d$ направлений.
%    \For 1 выполнять до $d-1$
%        \State Используя значения сетки и полученный на предыдущем шаге $\hat q_{i - 1}$, сгенерировать  $\hat q_i$.
%    \EndFor
%    \While{число итераций < $N_{TT}$}
%        \For 1 выполнять до $d-1$
%            \State На основе $\hat q_{i - 1}$ и $\hat q_i$ сгенерировать множество потенциальных решений $M$ и обновить сдвиг $\alpha$.
%            \State Запомнить лучшее решение $q_{best}$.
%            \State Представить массив значений функции $h(\tilde q)$, $\tilde q \in M$ в виде тензора.
%            \State Вычислить приближение тензора в ТТ-формате.
%            \State Используя значения сетки и $\hat q_{i - 1}$, сгенерировать $\hat q_i$.
%        \EndFor
%    \EndWhile
%\end{algorithmic}
%\end{algorithm}

The problem of minimizing the target functional~(\ref{eq:IP_func}) is planned to be solved by the Tensor Train (TT) global optimization method~\cite{Zheltkova_2018}, the algorithm of which is presented below.

\begin{algorithm}
\renewcommand{\thealgorithm}{\hspace{-5pt}}
\caption{of TT method}
\begin{algorithmic}[1]
    \Require{Lower and upper bounds of the solution space $b_{min}$ and $b_{max}$, number of parameters (dimensionality of the solution space) $d$, number of nodes along all directions $n$, maximum possible rank of wagons $r_{max}$, number of iterations $N_{TT}$, initial shift of the functional $\alpha$, mapping function $h(J(q) - \alpha)$.}
    \State Introduce a mesh with $n$ nodes along each of the $d$ directions.
    \For {1 to $d-1$}
        \State Using the grid values and the $\hat q_{i - 1}$ obtained in the previous step, generate $\hat q_i$.
    \EndFor
    \While{number of iterations < $N_{TT}$}
        \For {1 to $d-1$}
            \State Based on $\hat q_{i - 1}$ and $\hat q_i$, generate a set of potential solutions $M$ and update the shift $\alpha$.
            \State Remember the best solution $q_{best}$.
            \State Represent the array of values of the function $h(\tilde q)$, $\tilde q \in M$ as a tensor.
            \State Compute the tensor approximation in TT format.
            \State Using the grid values and $\hat q_{i - 1}$, generate $\hat q_i$.
        \EndFor
    \EndWhile
\end{algorithmic}
\end{algorithm}

\section*{CONCLUSION}

%Исследован вопрос идентифицируемости пространственной математической модели распространения быстротекущих эпидемий, основанной на законе действующих масс и диффузионных процессах. Прямая задача состоит в моделировании процесса распространения COVID-19 и определения функций плотностей восприимчивых, бессимптомных инфицированных, больных COVID-19, выздоровевших, госпитализированных, критических случаев, требующих подключение аппарата искусственной вентиляции легких и умерших в результате COVID-19. Построены алгоритмы численного решения прямой задачи, основанные на методе конечных элементов и методе конечных разностей. Свормклирована обратная задача, которая состоит в восстановлении вектора параметров, характеризующих процесс распространения COVID-19, по дополнительной информации о количестве выявленных, критических и умерших случаев за некоторый период времени. Алгоритм исследования вопроса идентифицируемости модели основан на глобальных методах анализа чувствительности Соболя и байесовском подходе, которые в совокупности позволили уменьшить границы вариации неизвестных параметров для дальнейшего решения обратной задачи. Показано, что для идентификации диффузионных коэффициентов, отвечающих за скорость перемещения индивидуумов в пространстве, необходимо использовать дополнительную информацию о процессе. А также предложен алгоритм решения обратной задачи, основанный на минимизации квадратичного целевого функционала методом глобальной оптимизации тензорного поезда.

The issue of identifiability of a spatial mathematical model of the spread of fast-moving epidemics based on the law of acting masses and diffusion processes is investigated. The direct problem consists of modelling the COVID-19 spreading process and determining density functions of susceptible, asymptomatic infected, COVID-19 patients, recovered, hospitalized, critical cases requiring artificial lung ventilation and COVID-19 deaths. Algorithms for numerical solution of the direct problem based on finite element method and finite difference method were constructed. The inverse problem, which consists of recovering the vector of parameters characterizing the process of COVID-19 spreading by additional information about the number of identified, critical and dead cases for some period of time, is developed. The algorithm for investigating the issue of model identifiability is based on global methods of Sobol sensitivity analysis and Bayesian approach, which together allowed to reduce the variation boundaries of unknown parameters for further solution of the inverse problem. It is shown that for identification of diffusion coefficients responsible for the rate of movement of individuals in space, it is necessary to use additional information about the process. And also an algorithm for solving the inverse problem based on minimization of the quadratic target functional by the global tensor train optimization method is proposed.


\begin{thebibliography}{99}

\bibitem{Kermack_McKendrick_1927}
\textit{Kermack~W.~O., McKendrick~A.~G.}
A contribution to the mathematical theory of epidemics~//
Proc. Roy. Soc. Lond. A. 1927. V.~115. P.~700--721.

\bibitem{KOI_KSI_2024}
\textit{Krivorotko~O., Kabanikhin~S.}
Artificial intelligence for COVID-19 spread modeling~//
J. Inverse Ill-Posed Probl. 2024. V.~32, No.~2. P.~297--332.

\bibitem{KOI_KSI_2020}
\textit{Krivorotko~O.~I., Kabanikhin~S.~I., Zyatkov~N.~Yu., Prikhodko~A.~Yu., Prokhoshin~N.~M., Shishlenin~M.~A.}
Mathematical Modeling and Forecasting of COVID-19 in Moscow and Novosibirsk Region~//
Num. Anal. Appl. 2020. V.~23, No.~4. P.~332--348.

\bibitem{KolmPetrPiskunov_1937}
\textit{Kolmogorov~A.~N., Petrovskii~I.~G.,  Piskunov~N.~S.}
A study of the diffusion equation with increase in the
amount of substance, and its application to a biological problem~//
Byull. Mosk. Gos. Univ. Mat. Mekh. 1937. V.~1. No.~6. P.~1--26.

\bibitem{Murray_2007}
\textit{Murray~J.~D.} 
Mathematical Biology. N.~Y.: Springer-Verlag, 2007.

\bibitem{Viguerie_2020}
\textit{Viguerie~A., Veneziani~A., Lorenzo~G., Baroli~D., Aretz-Nellesen~N., Patton~A., Yankeelov~T., Reali~A., Hughes~T., Auricchio~F.} 
Diffusion-reaction compartmental models formulated in a continuum mechanics framework: application to COVID-19, mathematical analysis, and numerical study~//
Comput. Mech. 2020. V.~66. P.~1131--1152.

\bibitem{Aristov_2021}
\textit{Aristov~V.~V., Stroganov~A.~V., Yastrebov~A.~D.} 
Simulation of spatial spread of the COVID-19 pandemic on the basis of the kinetic-advection model~//
Physics. 2021. V.~3, No~1. P.~85--102.

\bibitem{Barwolff_2021}
\textit{B\"{a}rwolff~G.} 
A local and time resolution of the COVID-19 propagation~--- a two-dimensional approach for Germany including diffusion phenomena to describe the spatial spread of the COVID-19 pandemic~//
Physics. 2021. V.~3, No.~3. P.~536--548.

\bibitem{Lau_2022}
\textit{Lau~Z., Griffiths~I.~M., English~A., Kaouri~K.} 
Predicting the spatially varying infection risk in indoor spaces using an efficient airborne transmission model~//
Proc. R. Soc. A. 2022. V.~478, No.~2259. Article~20210383.

\bibitem{Boscheri_2021}
\textit{Boscheri~W., Dimarco~G., Pareschi~L.}
Modeling and simulating the spatial spread of an epidemicthrough multiscale kinetic transport equations~//
Math. Models Methods Appl. Sci. 2021. V.~31, No.~6. P.~1059--1097.

\bibitem{Takacs_2021}
\textit{Tak\'acs~B.~M., Farag\'o~I., Horv\'ath~R., Repov\u{s}~D.}
Qualitative properties of space-dependent SIR models with constant delay and their numerical solutions~//
arXiv. 2022. arXiv:2112.06808.

\bibitem{KOI_ZNY_Eurasian_2022}
\textit{Krivorotko~O., Zyatkov~N.} 
Data-driven regularization of inverse problem for SEIR-HCD model of COVID-19 propagation in Novosibirsk region~//
Eurasian J. Math. Comput. Appl. 2022. V.~10, No.~1. P.~51--68.

\bibitem{covid_data}
\url{https://ai-biolab.ru/data}

\bibitem{Kabanikhin_2009}
\textit{Kabanikhin~S.} 
Definitions and examples of inverse and ill-posed problems~//
J. Inverse Ill-Posed Probl. 2009. V.~16, No.~4. P.~317--357.

\bibitem{Cukier} 
\textit{Cukier~R.~I., Fortuin~C.~M., Schuler~K.~E., Petschek~A.~G., Schaibly~J.~H.} 
Study of the sensitivity of coupled reaction systems to uncertainties in rate coefficients~// 
J. Chem. Phys. 1973. V.~59. P.~3873--3878.

\bibitem{Sobol} 
\textit{Sobol~I.~M.} 
Sensitivity analysis for non-linear mathematical models~// 
Math. Model. Comput. Exp. 1993. V.~1. P.~407--414.

\bibitem{Saltelli_2010} 
\textit{Saltelli~A., Annoni~P., Azzini~I., Campolongo~F., Ratto~M., Tarantola~S.} 
Variance based sensitivity analysis of model output. Design and estimator for the total sensitivity index~// 
Comput. Phys. Commun. 2010. V.~180. P.~259--270.

\bibitem{emulator} 
\textit{Andrianakis~I., Vernon~I.~R., McCreesh~N., McKinley~T.~J., Oakley~J.~E., Nsubuga~R.~N., Goldstein~M., White~R.~G.} 
Bayesian history matching of complex infectious disease models using emulation: a tutorial and a case study on HIV in Uganda~// 
PLoS Comput. Biol. 2015. V.~11. No.~1. Article~e1003968.

\bibitem{LHD} 
\textit{McKay~M.~D., Beckman~R.~J., Conover~W.~J.} 
A сomparison of three methods for selecting values of input variables in the analysis of output from a computer code~// Technometrics. 1979. V.~21. P.~239--254. 

\bibitem{lfbgs} 
\textit{Malouf~R.} 
A comparison of algorithms for maximum entropy parameter estimation~// 
Proc. Sixth Conf. Natur. Lang. Learn. 2002. P.~49--55.

\bibitem{Pukelsheim}  
\textit{Pukelsheim~F.} 
The three sigma rule~// 
Am. Stat. 1994. V.~48. P.~88--91.

\bibitem{Zheltkova_2018}
\textit{Zheltkova~V.~V., Zheltkov~D.~A., Grossman~Z., Bocharov~G.~A., Tyrtyshnikov~E.~E.}
Tensor based approach to the numerical treatment of the parameter estimation problems in mathematical immunology~//
J. Inverse Ill-Posed Probl. 2018. V.~26. No.~1, P.~51--66.

\end{thebibliography}
\end{document}